\documentclass[12pt]{amsart}
\usepackage[T2A]{fontenc}
\usepackage[cp1251]{inputenc}
\usepackage[english]{babel}
\usepackage{amsmath,latexsym,amsthm,amsfonts,tipa,upgreek}
\usepackage{amssymb,amscd,graphpap, stmaryrd}

\begin{document}

\title{Some more algebra on ultrafilters in metric spaces}
\author{Igor Protasov}
\keywords{Large, small, thin, thick, sparse and scattered subsets of groups; descriptive complexity; Boolean algebra of subsets of a group;  Stone-$\check{C}$ech compactification; ultracompanion; Ramsey-product  subset of a group; recurrence; combinatorial derivation.}

\maketitle

\begin{abstract}
We continue  algebraization of the set of ultrafilters on a metric spaces initiated in [6].
In particular, we define and study metric counterparts of prime, strongly  prime and right cancellable ultrafilters from the Stone-$\check{C}$ech compactification of a discrete group as a right topological semigroup [3]. Our approach is based on the concept of parallelity  introduced in the context of balleans in [4].
\vspace{3 mm}

{\bf 2010 MSC}: 22A15, 54D35, 54E15
\vspace{3 mm}

{\bf Keywords} : metric  space, ultrafilter,  ball invariance, parallelity, prime and strongly prime  ultrafilters.
\end{abstract}

\section{Introduction}

Let $X$ be a discrete space, and let $\beta X$ be the Stone-$\check{C}$ech compactification of $X$. We  take the point of $\beta X$ to be the ultrafilters on $X$, with the point of $X$ identified
with principal ultrafilters, so $X^{\ast} =\beta X \setminus X$ is the set of
all free ultrafilters. The topology of $\beta X$ can be defined by
stating that the sets of the form
$\bar{A}=\{p\in \beta X: A\in  p\}$,
where $A$ is a subset of $X$, are base for the open sets.
Given a filter $\varphi$ on $X$, the set $\bar{\varphi} = \bigcap\{\bar{A} : A\in \varphi \}$
is closed in $\beta X$, and every non-empty closed subset of $\beta X$  can be obtained in this way.

If $S$ is a discrete semigroup, the semigroup multiplication has a natural extension to $\beta S$, see
 [3, Chapter 4]. The compact right topological semigroup has very rich algebraic structure and a plenty
 applications to combinatorics, topological algebra and
 functional analysis, see [1, 2, 3, 7, 8, 10]. To get the
 product $pq$  of  $p,q \in \beta S$,  one can take an
 arbitrary $P\in p$, and for each $x\in P$, pick $Q_{x}\in q$. The
 $\bigcup_{x\in P} xQ_{x}\in pq$ and these subsets form a basis of $pq$.

In [6], given a metric space $X$, we endowed $X$
with the discrete topology, introduced and characterized the metric counterparts in $\beta X$ of minimal
left  ideals and the closure of the minimal ideal in   $\beta S$.

In this note, we continue  algebraization of
$\beta X$, define and describe  the metric analogues of prime, strongly prime and right cancellable ultrafilters from $\beta G$, $G$  is a discrete  group. We
recall that  an ultrafilter $p\in G^{\ast}$ is {\it prime} if
 $p \notin G^{\ast}G^{\ast}$, and $p$  is strongly  {\it prime}  if $p\notin cl G^{\ast}G^{\ast}$.
An  ultrafilter  $p\in G^{\ast}$  is called {\it right cancellable} if, for any $q, r\in \beta G$, $qp=rp$ implies $q=r$.

The key observation: to detect whether $p\in G^{\ast}$
is  prime or strongly prime, we do not need to know how to multiply any two individual ultrafilters
but only what is the set $G^{\ast} q$, $q\in G^{\ast}$. Indeed,
$p$ is prime if and only if $p\notin G^{\ast} q$ for each $q\in G^{\ast}$.
If $G$ is countable then $p\in G^{\ast}$ is right cancellable
if and only if  $p\notin G^{\ast} p$, see [3, Theorem 8.18].
But the natural metric counterpart of $G^{\ast} p$  in $\beta X$
can be defined by means of the parallelity  relation  on ultrafilters introduced in  [4] for the
general case of balleans, and applied for algebraization
of $\beta X$, $X$ is a metric space, in [6].

\section{Ball invariance and parallelity}

Let $(X,d)$  be a metric space. For any $x\in X$,  $A\subseteq  X$,  $r\in \mathbb{R}^{+}$,  $\mathbb{R}^{+} = \{r\in \mathbb{R}:r \geq 0\}$,  we denote $$B(x,r)= \{y\in X: d(x, y)\leq r\},  \   B(A,r)= \bigcup_{a\in A} B(a,r).$$

Given an ideal $\mathcal{I}$  in the Boolean algebra $\mathcal{P} _{X}$  of all subsets of $X$,  ($A, B\in \mathcal{I} $, $C\subseteq  A\Longrightarrow A\bigcup B\in \mathcal{I}, \ C\in  \mathcal{I}$), we say that $\mathcal{I}$
is {\it ball invariant} if, for every $A\in \mathcal{I}$  and
 $r\in \mathbb{R}^{+} $, we have   $B(A, r)\in \mathcal{I}$.  If $\mathcal{I}$  is ball invariant and $\mathcal{I}\neq \{\emptyset\}$  then $\mathcal{I}$ contains the ideal $\mathcal{I} _{b}$  of all bounded subsets of $X$. A
 subset $A$  of $X$ is called {\it bounded} if $A\subseteq  B(x,r)$  for some $x\in X$  and $r\in \mathbb{R}^{+}$.

We say that a filter $\varphi$  on $X$ is {\it  ball invariant} if,
for every $A\in \varphi$ and $r\in \mathbb{R}^{+} $,  there exists  $C\in \varphi$ such that $B(C, r)\subseteq A$.

An ideal $\mathcal{I}$ is called   proper if $\mathcal{I}\neq \mathcal{P}_{X}$. For a proper ideal in $\mathcal{P}_{X}$, we denote by $\varphi_{\mathcal{I}}$ the filter  $\{X\setminus  A:  A \in  \mathcal{I}\}$ and put $A^{\wedge} = \bar{\varphi}_{\mathcal{I}}$ so $$A^{\wedge} =\{p\in \beta X: X\setminus A\in p\}.$$

We remind the reader that $X$ in $\beta X$ is endowed with the discrete topology
and use the parallelity equivalence on $\beta X$ defined in [6] by the rule: $p || q$ if and only if
there  exists
$r\in \mathbb{R}^{+}$  such that $B(P, r)\in q$ for each $P\in p$.
A subset $S$ of $\beta X$ is called {\it invariant} (with respect to the parallelity equivalence) if, for all $p,q \in  \beta X$,  $p\in S$
 and $p || q$ imply $q\in S$.

\vskip 7pt

{\bf Proposition 1}. {\it For a  proper ideal $\mathcal{I}$ in $\beta X$, the following statements  are equivalent:

$(i)$  $\mathcal{I}$ is ball invariant;

$(ii)$  $\varphi _{\mathcal{I}}$  is ball invariant;

$(iii)$  $\mathcal{I}^{\wedge}$ is invariant.}
\vskip 7pt

{\it Proof.}  The equivalence $(i) \Longleftrightarrow (ii)$  is evident. To prove $(ii)\Longrightarrow  (iii)$,   let $p\in \mathcal{I}^{\wedge}$  and $q || p$.
We choose $r \geq  0$  such that $B(P, r)\in q$   for each $P\in p$. Given an arbitrary $Y\in \varphi_{\mathcal{I}}$,
we
choose $Z\in  \varphi_{\mathcal{I}}$  such that $B(Z, r) \subseteq Y$.  Then $Z\in p$  and $B(Z, r)\in q$  so $Y\in q$ and $q\in \mathcal{I}^{\wedge}$.

To see that $(iii)\Longrightarrow  (i)$, we assume the contrary and choose
$Y\in \varphi_{\mathcal{I}}$  and $r\geq  0$  such that $B(Z, r)\setminus  Y\neq \emptyset$  for each $r \geq 0$.
Then we take $q\in \beta X$ such that $B(Z, r)\setminus Y \in q$  for each $Z\in \varphi_{\mathcal{I}}$.
   By [6, Lemma 2.1], there exists $p\in \varphi_{\mathcal{I}}$  such that
   $q \mid\mid p$.  Since $q  \notin  \mathcal{I}^{\wedge} $, we get a contradiction.  $ \ \ \Box$
   \vskip 7pt

In what follows, we suppose that every metric space $X$ under  consideration is unbounded, put
$$   X^{\sharp} = \{ p\in \beta X: \ every \ member \ P\in p  \ is \  unbounded \ in  \  X\}$$
and note that $X^{\sharp} $ is a closed invariant subset of $\beta X$.
\vskip 7pt

We say that a subset $A$ of $X$ is

 \begin{itemize}
\item{}   {\it large}  if $X= B(A, r)$  for some $r\geq 0$;

\item{}   {\it small}  if  $L \setminus A$  is large for every large subset $L$;

\item{}   {\it thick}  if, for every  $r \geq 0$, there exists $a\in A$  such  that $B(a,r)\subseteq A$;

\item{}   {\it prethick}  if  $B(A, r)$ is thick  for some  $r\geq 0$.
\end{itemize}
\vskip 10pt

The family $Sm _{X}$ of all small subsets of $X$ is an ideal  in $\mathcal{P}_{X}$, and a subset $A$  is small if and only if $A$  is not  prethick [7, Theorems 11.1 and 11.2].
\vskip 7pt

{\bf Proposition 2}. {\it For every metric  space $X$, the ideal $Sm _{X}$ is ball invariant and

$ Sm _{X}^{\wedge}=cl \{ \bigcup \{ K: \ K $   is  a  minimal  non-empty  closed  invariant subset  of } $X^{\sharp} \}\}$

\vskip 7pt

{\it Proof.}  The second statement is the dual form of Theorem  3.2 from [6].  Assume that $A$ is small  but $B(A, r)$ is not small for some $r\geq 0$.  Then   $B(A, r)$  is prethick so there is
$m\geq 0$  such  that $B(B(A,r), m)$  is thick. It follows that $A$ is prethick and we get a contradiction.  $ \  \Box$

For every   metric space $X$, by [6, Corollary 3.1], the set of all minimal non-empty closed invariant subset of $X$
 has cardinality $2 ^{2 ^{asden X}}$, where $ asden X= min \{ |Y| : Y$ is a large subset of $X\}$.   Applying Proposition 1, we get $2 ^{2 ^{asden X}}$  maximal proper ball invariant ideals in $\mathcal{P}_{X}$.

\vskip 7pt

{\bf Proposition 3}. {\it  Let $\mathcal{I}$  be a ball invariant ideal in  $\mathcal{P}_{X}$ such that
$\mathcal{I}\neq \mathcal{I}_{b}$,   $ \ \mathcal{I}_{b}$  is the ideal of all bounded subsets of $X$.  Then there exists a ball invariant ideal  $\mathcal{J}$ such that $\mathcal{I}_{b}\subset \mathcal{J}\subset  \mathcal{I}$}.
\vskip 7pt

 {\it Proof.}   We  take an unbounded subset      $A\in  \mathcal{I}$  and choose a sequence $(a_{n})_{n\in \omega}$  in $A$  such that $B(a_{n}, n)\bigcap B(a_{m}, m) = \emptyset$
 for all distinct $n, m \in \omega$. We put  $A_{0} = \{a_{2n}: n\in\omega\}$, $A_{1} = \{a_{2n+1}: n\in\omega\}$
  and   denote by $\mathcal{J}$  the smallest  ball invariant ideal such
  that $A_{0} \in \mathcal{I}$. Then $Y\in \mathcal{I}$ if and only if $Z\subseteq B(Y_{0}, m)$ for
  some $m\in \omega$. By the choice  of  $(a_{n})_{n\in\omega}$, $Y\setminus B(Y_{0}, m)\neq\emptyset$
  for each $m\in \omega$,  so  $Y_{1} \notin  \mathcal{J}$  but  $Y_{1}\in \mathcal{I}$.  $ \ \ \Box$

\section{Prime and strongly prime ultrafilters}

For each $q\in X^{\sharp}$, we denote $q^{=}= \{r\in X^{\sharp}: r || q\}$ and
 say that $p\in X^{\sharp}$ is {\it divisible} if there exists $q\in X^{\sharp}$ such
 that $\bar{P}\bigcap q^{=}$ is  infinite for each $P\in p$.  An ultrafilter
 $p\in X^{\sharp}$ is called {\it prime} if $p$ is not divisible, and
 {\it strongly prime} if $p$ is not in the closure of the set of all divisible    ultrafilters.

A subset $A$ of $X$ is called {\it sparse } if $\bar{A} \bigcap q^{=}$  is finite for  each $q\in X^{\sharp}$.
We denote by $Sp _{X}$ the family of all
sparse subsets of $X$  and observe that $Sp _{X}$  is an ideal in $\mathcal{P}_{X}$.
\vskip 7pt

{\bf Proposition 4}. {\it An ultrafilter $p\in X^{\sharp}$ is strongly prime if and only if there exists $A\in Sp_{X}$  such that $A\in p$  so $Sp_{X}^{\wedge}= cl \ \mathcal{D}$, where $\mathcal{D}$  is the set of all divisible ultrafilters.}
\vskip 7pt

{\it Proof.} Assume that each member $P\in p$ is not sparse and choose $q\in X^{\sharp}$ such that $\bar{P}\bigcap q^{=}$ is infinite. We
 take an arbitrary limit point $r$  of the set $\bar{P}\bigcap q^{=}$ . Then $P\in r$ and $r$ is divisible so $p\in cl \mathcal{D}$  and $p$  is not strongly prime.

On the other hand, if $A$  is sparse and $A\in p$  then $\bar{A} \bigcap \mathcal{D}=\emptyset$ and $p\notin cl \mathcal{D}$.
 $ \ \ \Box$
\vskip 7pt

A subset $A$ of $X$ is called {\it thin} if, for every $r \geq  0$, there exists a bounded subset $V$ of $X$  such that $B(a,r)\bigcap  A = \{a\}$  for each $a\in A\backslash V$.
\vskip 7pt

{\bf Proposition 5}. {\it If  $p\in X^{\sharp}$  and some member $P\in p$  is thin
then $p$ is strongly prime.
\vskip 7pt

Proof.}  By [6, Theorem 4.3],  $P$ is thin if and only if  $|\bar{P}\bigcap q^{=}|\leq 1$  for each $q\in X^{\sharp}$
 so we can apply Proposition 4.
  $ \ \ \Box$
 \vskip 7pt

Since every unbounded subset of $X$  contains some unbounded thin subset, we conclude that the set of all strongly prime ultrafilters is dense in $X^{\sharp}$.

Is the ideal $Sp _{X}$  ball invariant? In Proposition 7, we give a negative example. In Proposition 6, we describe a class of metric spaces for which the answer is positive.

A metric space $X$ is called {\it uniformly locally finite} if, for every $r\geq 0$,  there exists $m\in \mathbb{N}$  such that $|B(x,r)|\leq m$ for each $x\in X.$
\vskip 7pt

{\bf Proposition 6}. {\it If a metric space $X$ is uniformly locally finite then the ideal  $Sp _{X}$ is ball invariant.
\vskip 7pt

Proof.} By [5, Theorem 1],  there exists a countable group $G$ of permutations of $X$  such  that
\vskip 7pt

$(1)$ for each $r\geq 0$, there exists a finite subset $F$ of $G$ such that $B(x,r)\subseteq F(x)$  for each $x \in X$, where $F(x)=\{ g(x): g\in F\}$;
\vskip 7pt

$(2)$ for every finite subset $F$ of  $G$, there exists $r\geq 0$ such that  $F(x)\subseteq B(x,r)$  for each $x\in X$.
\vskip 10pt

It follows that, for $p,q\in  X^{\sharp}$, $p|| q$ if and only if
there exists $g\in G$ such that $q= g(p)$, where  $g(p)= \{ g(P): P\in p\}$.

Now let $A$ be a sparse subset of $X$ and $r\geq 0$.  We
choose $F$ satisfying $(2)$ so $B(A, r)\subseteq  F(A)$, where
$F(A)=\bigcup _{g\in F} g(A)$. We take an arbitrary $q\in X^{\sharp}$.
Since $A$ is sparse, $q^{=}\bigcap \bar{A}$  is finite.
Then $q^{=}\bigcap \overline{B(A, r)}\subseteq \bigcup_{g\in F}(q^{=}\bigcap  \overline{g(A)})$.
Since $|q^{=}\bigcap \overline{g(A)}|$
 $=|(g^{-1}q )^{=}\bigcap A|$
 and $A$ is sparse, $|q^{=}\bigcap \overline{B(A,r)}|$
 is finite and $B(A, r)$  is sparse.   $ \  \Box$
\vskip 7pt

{\bf Proposition 7}. {\it Let $\mathbb{Q}$  be the set of rational numbers endowed with the metric $d(x,y)= |x-y|$.  The ideal $Sp _{\mathbb{Q}}$ is not ball invariant.
\vskip 7pt

Proof.} We put $A=\{ 2^{n}: n \in \mathbb{N} \}$. By Proposition 5, $A$ is sparse.
 We take an arbitrary free ultrafilter $q\in \bar{A}$.  Then $B(A,1)\in x+q$ for each $x\in [0,1]$. Since $x+q ||q$,
 $ q^{=} \bigcap  \overline{B(A,1)}$ is infinite so $B(A,1) $  is not sparse. $ \  \Box$

We say that an ultrafilter $p\in X^{\sharp}$  is {\it discrete} if each $q\in p^{=}$ is an isolated point in the set $p^{=}$.  In view
of [3, Theorem 8.18], a discrete ultrafilter can be considered as a counterpart of a right cancellable ultrafilter. Clearly, if each $q\in p^{=}$ is prime then $p$ is discrete.
\vskip 7pt

{\bf Proposition 8}. {\it There exist two ultrafilters $p,q\in \mathbb{Q}^{\sharp}$  such that $p|| q$,  $p$ is isolated in $p^{=}$  but $q$  is not isolated in $p^{=}$.
\vskip 7pt

Proof.}  For each $n\in \mathbb{N}$, we  put  $A_{n}= \bigcup_{m\geq n}[2^{m, 2^{m}+1}]$
 and take a maximal filter $q$  such that $A_{n}\in q$,  $n\in \mathbb{N}$ and each member $A\in q$ is somewhere dense,
i.e. the closure of $A$ in $\mathbb{Q}$ has non-empty interior.
 It is easy to see that $q$ is an ultrafilter and $q$  has a basis consisting of subsets without isolated points. We consider the mapping
$f: \mathbb{Q }\longrightarrow \mathbb{Q}$ defined by $f(x)=  \lfloor x\rfloor $,  where $\lfloor x\rfloor $ is the nearest from the left integer to $x$.  The set $\{f(\mathcal{U}): \mathcal{U}\in q\}$
is a basis for some uniquely determined ultrafilter $p$ such that $\{2^{n} : n \in \mathbb{N} \} \in p$. Clearly $p||q$ and, by Proposition 5, $p$ is isolated in $p^{=}$.

We   show that $q$ is not isolated in $p^{=}$.  We take
an arbitrary  $Q\in q$ such that  $Q$ has no isolated points, $f(Q)\subseteq \{2^{n}: n\in \mathbb{N}\} $ and
choose an arbitrary mapping $h: f(Q) \longrightarrow Q$ such
that $h(2^{n})\in  [2^{n}, 2^{n} +1]$ for each $2^{n}\in f(Q)$. We
denote by $q_{h}$ the ultrafilter with the basis
$\{h(\mathcal{V}): \mathcal{V}\in p\}$. Then $Q\in q_{h}$ and $q_{h} \mid\mid p$. Since
$Q$  has no isolated points, we have countably many
ways to choose $h$  and get countably many distinct
ultrafilters from $p^{=}\bigcap \bar{Q}$.   $ \  \Box$

A subset $A$ of $X$ is called {\it disparse} if $\bar{A }\bigcap p^{=}$
is discrete for each $p\in X^{\sharp}$. The family $dSp _{X}$ of all
disparse subsets of $X$ is an ideal in $\mathcal{P}_{X}$  and we get
the following evident
\vskip 7pt

{\bf Proposition 9}. {\it   For every metric space $X$,  $dSp_{X}^{\wedge}$ is the
of all ultrafilters $p\in X^{\sharp}$ such that  $p^{=}$ has no isolated points.}

\vskip 7pt

{\bf Proposition 10}. {\it  For every $p\in X^{\sharp}$, the set $p^{=}$  is
nowhere dense in $X^{\sharp}$.
\vskip 7pt

Proof.} We take an arbitrary $A\in p$ and coming back to the proof of Proposition  3, consider the subsets $A_{0}$, $A_{1}$ of $A$. If $A_{0} \in q$, $A_{1}\in r$  then $q$ and $r$ are not
parallel. Then either $\overline{A_{0}}\bigcap p^{=}$ or $\bar{A_{1}}\bigcap p^{=}=\emptyset$. $ \  \Box$

\section{Ballean context}

Following [7, 8], we say that a {\it ball structure}
is a triple $\mathcal{B}=(X, P, B)$, where $X, P$ are non-empty sets
and, for every $x\in X$ and $\alpha\in P$,  $B(x, \alpha)$ is a subset of
$X$  which is called a {\it ball of radius} $\alpha$  around $x$.  It is
supposed that $x\in B(x, \alpha)$ for all $x\in X$ and $\alpha\in P$. The
set $X$ is called the {\it support} of $\mathcal{B}$, $P$ is called the {\it set of radii.}

Given any $x\in X$, $A \subseteq X$, $\alpha\in P$, we set
$$B^{\ast}(x,\alpha)=\{y\in X: x\in B(y,\alpha)\},  \   B(A,\alpha)=\bigcup _{a\in A}B(a,\alpha).$$

A ball structure $\mathcal{B}=(X P, B)$  is called a {\it ballean} if

\begin{itemize}

\item{} for any $\alpha,\beta\in P$, there exist $\alpha^{\prime},\beta^{\prime}\in P$ such that,
for every $x\in X$,
$$B(x,\alpha)\subseteq B^{\ast}(x,\alpha^{\prime}), \ \  B^{\ast}(x,\beta)\subseteq B(x,\beta^{\prime});$$

\item{} for any $\alpha,\beta\in P$, there exists $\gamma\in P$ such that,
for every $x\in X$,
$$ B(B(x,\alpha),\beta)\subseteq B(x,\gamma);$$

\item{} for any $x,y\in X$, there is $\alpha\in P$ such that   $y\in B(x,\alpha)$.
\end{itemize}
\vskip 10pt

A ballean $\mathcal{B}$ on $X$ can also be determined in terms
of entourages of the diagonal $\triangle_{X}$ in $X\times  X$  (in this case
it is called a {\it coarse structure}  [9]) and can be
considered as an asymptotic  counterpart of a uniform topological space.

Every metric space $(X, d)$ defines the  ballean
$(X, \mathbb{R}^{+} , B_{d})$, where $B_{d}(x,r)=
\{ y\in X: d(x,y)\leq r\}$. For
criterion of metrizability of balleans see [8 , Theorem 2.1.1].

We observe that all definitions in this paper do not use the metric on $X$ directly but only balls so
can be literally rewritten for  any ballean in place
of metric space. Moreover, a routine verification
ensures that Propositions 1, 2, 4, 5, 9 remain true for any balleans.

Let $G$ be  a group with the identity $e$. We
denote by $\mathcal{F}_{G}$  the family of all finite subsets of $G$  containing $e$  and get the group ballean  $\mathcal{B}(G)=(G, \mathcal{F}_{G}, B) $, where $B(g, F)= Fg$  for all $g\in G$,  $F\in \mathcal{F}_{G}$.
We note that $G^{\sharp}=G^{\ast}$  and, for $p, q\in G^{\ast}$, $p \mid\mid q$
if and only if $q = gp$  for some $g\in G$.  Hence, $p^{=}=Gp$,
$cl p^{=}= (\beta G) p$
and the minimal non-empty  closed
invariant subsets   in $G^{\sharp}$  are precisely the minimal
left ideals of the semigroup  $\beta G$.   The  ballean and
semigroup notions of divisible, prime and strongly prime ultrafilters coincide.

\vspace{3 mm}
CONTACT INFORMATION

I.~Protasov: \\
Faculty of Computer Science and Cybernetics  \\
        Kyiv University  \\
         Academic Glushkov pr. 4d  \\
         03680 Kyiv, Ukraine \\ i.v.protasov@gmail.com

\end{document}